%
%
%
%

\def\cent#1#2{{\bf C}_{#1}(#2)}
\def\sbs{\subseteq}
\def\irr#1{{\rm Irr}(#1)}
\def\ref#1{{\bf[#1]}}
\def\iitem#1{\goodbreak\par\noindent{\bf #1}}
\def\qed{{~~~\vrule height .75em width .4em depth .3em}}
\def\zent#1{{\bf Z}(#1)}
\def\ker#1{{\rm ker}(#1)}
\def\hat{\widehat}

\magnification = \magstep1
{\nopagenumbers
\font\bbf = cmbx12
\font\teb=cmmib10 scaled \magstep2

\centerline{\bbf CAMINA {\teb p}-GROUPS THAT ARE}
\medskip
\centerline{\bbf GENERALIZED FROBENIUS COMPLEMENTS}
\bigskip\bigskip
\centerline{by}
\medskip

\centerline{\bf I. M. Isaacs}
\smallskip
\centerline{Mathematics Department}
\centerline{University of Wisconsin}
\centerline{480 Lincoln Dr.}
\centerline{Madison, WI ~~~53706}
\centerline{USA}
\smallskip
\centerline{E-mail: isaacs@math.wisc.edu}
\medskip
\centerline{and}
\medskip

\centerline{\bf Mark L. Lewis}
\smallskip
\centerline{Department of Mathematical Sciences}
\centerline{Kent State University}
\centerline{Kent, OH ~~~44242}
\centerline{USA}
\smallskip
\centerline{E-mail: lewis@math.kent.edu}

\vglue 1.5truein

\centerline{\bf ABSTRACT}
\medskip
{\narrower
Let $P$ be a Camina $p$-group that acts on a group $Q$ in such a
way that $\cent Px \sbs P'$ for all nonidentity elements $x \in Q$. We
show that $P$ must be isomorphic to the quaternion group $Q_8$. If
$P$ has class $2$, this is a known result, and this paper corrects a
previously published erroneous proof of the general case.\par}

\vfil\noindent

\eject}

Let $G$ be a finite group that is neither perfect nor abelian, and
recall that $G$ is said to be a {\bf Camina group} if every nontrivial
coset of $G'$ is a conjugacy class of $G$, or equivalently, every
nonlinear irreducible character of $G$ vanishes on $G - G'$. (That
these conditions really are equivalent is immediate from the fact that
for $x \in G$, we have $|\cent Gx| = \sum|\chi(x)|^2$, where the sum
runs over $\chi \in \irr G$.)

The purpose of this note is to prove the following, which appeared
with an incorrect proof as Theorem~2 of \ref3.
\medskip

\iitem{THEOREM.}~~{\sl Let $P$ be a Camina $p$-group, and
suppose that $P$ acts on a nontrivial group $Q$ in such a way that
$\cent Px \sbs P'$ for all nonidentity elements $x \in Q$. Then $P$ is
the quaternion group of order $8$, and the action is Frobenius.}
\medskip

This theorem was used in \ref3 to give an alternative proof of a key
step in the classification of Camina groups given in \ref1, where Dark
and Scoppola  proved that a Camina group must be either a
$p$-group or a Frobenius group whose complement is either cyclic
or $Q_8$. Unfortunately, as is explained in \ref3, the Dark-Scoppola
proof ultimately relies on a flawed argument in \ref1. Combining our
result with Theorem~1 of \ref3, we now have what we hope is a
correct (and simplified) proof of the Dark-Scoppola classification.

We begin with a fairly standard general result.
\medskip

\iitem{LEMMA.}~~{\sl Let $A$ and $B$ be finite abelian groups, and
suppose that there exists a nondegenerate bimultiplicative map
$f: A \times B \to C$, where $C$ is a finite cyclic group. (This means
that $f$ is a homomorphism in each variable and that the only
elements $a \in A$ and $b \in B$ such that $f(a,B) = 1$ or $f(A,b) = 1$
are the identities of $A$ and $B$.) Then $A \cong B$.}
\medskip

\iitem{Proof.}~~Let $\mu$ be a faithful linear character of $C$. For
each element $a \in A$, let $\lambda_a$ be the function on $B$
defined by $\lambda_a(b) = \mu(f(a,b))$. It is easy to check that the
map $a \mapsto \lambda_a$ is an injective homomorphism from $A$
into the group $\hat B$ of linear characters of $B$. Then
$|A| \le |\hat B| = |B|$, and by symmetry, we have $|A| = |B|$. It
follows that $A \cong \hat B \cong B$, as required.\qed
\medskip

Next, we present a few easy results about Camina $p$-groups.
Stronger versions of many of these are known, but they are scattered
over a number of papers. (For example, see \ref1 and the references
there.) Since the facts that we need can be established with
elementary arguments, it seems reasonable to present the proofs
here.

In the following, $P$ is always a Camina $p$-group.
\medskip

\iitem{PROPOSITION 1.}~~{\sl $P/P'$ is elementary abelian.}
\medskip

\iitem{Proof.}~~Let $x \in P$. If $x^p \not \in P'$, then $x \not \in P'$,
so $|P:\cent Px| = |P'| = |P:\cent P{x^p}|$, and thus
$|\cent Px| = |\cent P{x^p}|$. Since $\cent Px \sbs \cent P{x^p}$, we
deduce that $\cent Px = \cent P{x^p}$. Now let $z \in \zent P$ have
order $p$, and note that $z \in P'$, so $x$ and $zx$ are conjugate.
Then $x^t = zx$ for some element $t \in P$, and we have
$(x^p)^t = (x^t)^p = (zx)^p = x^p$. It follows that
$t \in \cent P{x^p} = \cent Px$, and this is a contradiction since
$x^t = zx \ne x$. Thus $x^p \in P'$, as required.\qed
\medskip

\iitem{PROPOSITION 2.}~~{\sl Suppose $P$ has nilpotence class
$2$. Then $|P:P'| > |P'|$.}
\medskip

\iitem{Proof.}~~Let $x \in P - P'$. Then $|P:\cent Px| = |P'|$, so
$|P:P'| = |\cent Px| > |P'|$, where the inequality is strict because
$x \in \cent Px$.\qed
\medskip

>From now on, we assume that $P$ has nilpotence class $3$. (We
mention that by the Theorem in \ref1, the nilpotence class of a
Camina $p$-group cannot exceed $3$. We will use this fact in the
proof of our main result.)
\medskip

\iitem{PROPOSITION 3.}~~{\sl $[P',P]$ is elementary abelian.}
\medskip

In fact, more is true: $P'$ is elementary abelian. (See the remarks
preceding and following Proposition~6, below.)
\medskip

\iitem{Proof of Proposition~3}~~Since $[P',P]$ is abelian and is
generated by elements of the form $[u,x]$, where $u \in P'$ and
$x \in P$, it suffices to show that $[u,x]^p = 1$. Since $[u,x]$ is
central in $P$, we have $[u,x]^p = [u,x^p] = 1$, where the final
equality holds since $P'$ is abelian, and by Proposition~1, it contains
both $u$ and $x^p$.\qed
\medskip

In the following, $Z = \zent P$ and $C = \cent P{P'}$, and we note
that $Z \sbs P' \sbs C$.
\medskip

\iitem{PROPOSITION 4.}~~{\sl Assume that $p > 2$. Then $C/Z$ is
elementary abelian.}
\medskip

Actually, it is not really necessary to assume that $p > 2$ here or in
Proposition~6. This is because by Theorem~3.1 of \ref4, Camina
$2$-groups never have nilpotence class exceeding $2$. (This fact too
will be used in our proof of the main result.) 
\medskip

\iitem{Proof of Proposition~4.}~~Let $c \in C$ and $x \in P$, and
write $c^x = cu$ and $u^x = uv$, where $u \in P'$ and
$v \in [P',P] \sbs Z$. For positive integers $n$, it follows by induction
that $c^{x^n} = cu^nv^{n(n-1)/2}$. Since $p > 2$ and $v^p = 1$ by
Proposition~3, we have $c^{x^p} = cu^p$. Also $x^p \in P'$ by
Proposition~1, and thus since $c \in C = \cent P{P'}$, we have
$c^{x^p} = c$, and thus $u^p = 1$. Now
$(c^p)^x = (c^x)^p = (cu)^p = c^p$, where the last equality holds
because $c$ and $u$ commute and $u^p = 1$. Since $x \in P$ was
arbitrary, it follows that $c^p \in Z$, and thus $C/Z$ has exponent
$p$.

To see that $C/Z$ is abelian, note that $[P,C,C] \sbs [P',C] = 1$, so it
follows by the three-subgroups lemma that $[C',P] = 1$, and thus
$C' \sbs Z$.\qed
\medskip

\iitem{PROPOSITION 5.}~~{\sl Assume that $Z$ is cyclic. Then
$|C:P'|$ is a square and $|C:P'| \ge p^2$.}
\medskip

\iitem{Proof.}~~Since $Z$ is cyclic, $P$ has a faithful irreducible
character $\chi$, and we argue that $\chi$ vanishes $P - Z$, and
thus $|P:Z| = \chi(1)^2$. To see this, observe that $\chi$ vanishes on
$P - P'$ since $P$ is a Camina group, so it suffices to show that
$\chi(x) = 0$ for $x \in P' - Z$. Since $x$ is noncentral, we can
choose $t \in P$ with $x^t = xz$, where $z \in Z$ is some nonidentity
element. Now $\chi_Z = \chi(1)\lambda$ for some faithful linear
character $\lambda$ of $Z$, and thus
$\chi(x) = \chi(x^t) = \chi(xz) = \lambda(z)\chi(x)$. Since
$\lambda(z) \ne 1$, it follows that $\chi(x) = 0$, as wanted.

Now commutation defines a nondegenerate bimultiplicative map from
$(P/C) \times (P'/Z)$ into the cyclic group $Z$, and thus
$|P:C| = |P':Z|$ by the lemma. Since $|P:Z|$ is a square, it follows
that $|C:P'|$ is a square. Also, $P/Z$ is a class $2$ Camina group, so
by Proposition~2, we have $|P:P'| > |P':Z| = |P:C|$, and thus
$C > P'$, and we have $|C:P'| \ge p^2$.\qed
\medskip

\iitem{PROPOSITION  6.}~~{\sl Assume that $Z$ is cyclic and
$p > 2$. Then $P'$ is elementary abelian.}
\medskip

As we have remarked, the assumption that $p > 2$ is redundant
here. In fact, as we will explain following the proof, the assumption
that $Z$ is cyclic is not really needed either.
\medskip

\iitem{Proof of Proposition 6.}~~By Proposition~5, we can choose an
element
$c \in C - P'$. Now let $u \in P'$ be arbitrary, and choose an element
$t \in P$ such that $c^t = uc$. Since $c^p \in Z$ by Proposition~4,
we have $c^p = (c^p)^t = (c^t)^p = (uc)^p = u^pc^p$, where the last
equality holds because $c$ and $u$ commute. Then $u^p = 1$, as
required.\qed
\medskip

To see why it is not really necessary to assume that $Z$ is cyclic in
Proposition~6, observe that it suffices to show for each character
$\chi \in \irr G$ that the group $(P/\ker\chi)'$ is elementary abelian.
Now $P/\ker\chi$ has a cyclic center, and it is either abelian, a
Camina group of class $2$ or a Camina group of class $3$. If it is
abelian, there is nothing to prove and if it has class $3$, its derived
subgroup is elementary by Proposition~6. Finally, is easy to see
using Proposition~1 that class $2$ Camina $p$-groups always have
elementary abelian derived subgroups.
\medskip

\iitem{Proof of Theorem.}~~We proceed by induction on $|P|$.
Observe first that the hypothesis guarantees that $P$ centralizes no
nonidentity element of $Q$, and thus $|Q| \equiv 1$ mod $p$, and
hence $Q$ is a $p'$-group. Since $P$ is a Camina $p$-group, its
nilpotence class is at least $2$. If $P$ has class $2$, the result
follows by Lemma~3.1 of \ref2, so we can assume that the class of
$P$ is at least $3$. Since Camina $2$-groups have nilpotence class
$2$ by Theorem~3.1 of \ref4, we deduce that $p > 2$, and we work
to obtain a contradiction. Also, by the Theorem in \ref1, Camina
$p$-groups have class at most $3$, and hence the class of $P$
must be $3$, exactly.

Since we can replace $Q$ by a nontrivial $P$-invariant subgroup, we
can assume that $Q$ has no proper nontrivial $P$-invariant
subgroup. It follows that $Q$ is an elementary abelian $q$-group for
some prime $q \ne p$, and thus we can view $Q$ as an irreducible
$F[P]$-module, where $F$ is the field of order $q$.

Let $K$ is the centralizer of the action of $G$ in the endomorphism
ring of $Q$. Then $K$ is a finite division ring, and hence by
Wedderburn's theorem, $K$ is a field. We can view $Q$ as a vector
space over $K$, and as such, $Q$ is an absolutely irreducible
$K[P]$-module. (We stress that we have not changed $Q$; what has
changed is our point of view.) Let $\chi$ be the irreducible
$q$-Brauer character of $P$ corresponding to the absolutely
irreducible $K[P]$-module $Q$, and observe that in fact,
$\chi \in \irr G$ since $q$ does not divide $|P|$.

We argue now that $\chi$ is faithful. Otherwise, there exists a
nontrivial central subgroup $U \sbs \ker\chi$, and hence $U$ acts
trivially on $Q$. Also, $P/U$ is nonabelian since $P$ has class $3$,
and it follows that $P/U$ is a Camina $p$-group, so we can apply the
inductive hypothesis to the action of $P/U$ on $Q$. Then $P/U$ has
order $8$, and this is a contradiction since we have established that
$p > 2$. 

Next, we show that no element of $P$ outside of $P'$ can have
order $p$. To see this, suppose that $H \sbs P$, where $|H| = p$ and
$H \cap P' = 1$. Since $P$ is a Camina group, $\chi$ vanishes on
the nonidentity elements of $H$, and hence $\chi_H$ has a principal
constituent. The action of $H$ on $Q$ is completely reducible,
however, and it follows that $H$ has nontrivial fixed points in $Q$,
and hence by hypothesis, $H \sbs P'$, and this is a contradiction.

Let $C = \cent G{P'}$ and $Z = \zent P$, and note that $Z$ is cyclic
since $\chi$ is faithful. Also $Z \sbs P'$, and thus $|Z| = p$ by
Proposition~6. If $c \in C - P'$, then $c^p \in Z$ by Proposition~4,
and since $c$ does not have order $p$, we deduce that $c^p$ is a
generator of $Z$.

Now $C/P'$ is elementary abelian by Proposition~1, and it has order
at least $p^2$ by Proposition~5, and thus we can choose two
elements $b,c \in C$ that generate distinct subgroups of order $p$
modulo $P'$. Then $b^p$ and $c^p$ are generators of $Z$, and we
can replace $c$ by a suitable power and assume that in fact,
$b^p = c^{-p}$. Now write $[b,c] = z$, so $z \in Z$ by Proposition~4.
Then $(bc)^p = b^pc^pz^{p(p-1)/2} = b^pc^p = 1$ because $p > 2$
and $z^p = 1$. Since $b$ and $c$ generate different subgroups of
order $p$ in $C/Z$, it follows that $bc \not\in P'$, and since $bc$
has order $p$, this is our final contradiction.\qed
\bigskip\bigskip
\noindent\frenchspacing
\centerline{\bf REFERENCES}
\bigskip

\parindent = 15pt
\noindent
\item{\bf1.} R. Dark and C. M. Scoppola, On Camina group of prime
power order, {\it J. Algebra} {\bf 181} (1996), 787--802.
\medskip

\noindent
\item{\bf2.} I. M. Isaacs, Coprime group actions fixing all nonlinear
irreducible characters, {\it Canad. J. Math.} {\bf 41} (1989), 68--82.
\medskip

\noindent
\item{\bf3.} M. L. Lewis, Classifying Camina groups: a theorem of
Dark and Scoppola, {\it Rocky Mountain J. Math.} {\bf 44} (2014),
591--597.
\medskip

\noindent
\item{\bf4.} I. D. Macdonald, More on $p$-groups of Frobenius type,
{\it Israel J. Math.} {\bf 56} (1986), 335--344.

\bye